\input amstex
\documentstyle{amsppt}
\magnification=1200
\hoffset=-0.5pc
\vsize=57.2truepc
\hsize=38truepc
\nologo
\spaceskip=.5em plus.25em minus.20em
\define\parti{d}
\define\ltime{\ltimes}
\define\barakont{1}
\define\batviltw{2}
\define\batvilfo{3}
\define\batavilk{4}

\define\canhawei{7}
\define\cheveile{8}
\define\gersthtw{9}
\define\geschthr{10}
\define\getzltwo{11}
\define\herzone{12}
\define\poiscoho{13}
\define\duality{14}
\define\bv{15}
\define\extensta{16}
\define\crosspro{17}
\define\twilled{18}
\define\liribi{19}
\define\huebstas{20}
\define\kosmathr{21}
\define\kosmafou{22}
\define\kosmafiv{23}
\define\kosmagtw{24}
\define\koszulon{25}
\define\liazutwo{26}
\define\luweinst{27}
\define\mackfift{28}
\define\mackxu{29}
\define\maclaboo{30}
\define\majidtwo{31}
\define\maninfiv{32}
\define\mokritwo{33}
\define\palaione{34}
\define\rinehone{35}
\define\schecone{36}
\define\stachone{37}
\define\stashnin{38}
\define\tianone{39}
\define\todortwo{40}
\define\wittetwe{41}
\define\xuone{42}
\define\zakrztwo{43}
\define\Bobb{\Bbb}
\define\fra{\frak}
\topmatter
\title
Differential Batalin-Vilkovisky algebras\\
arising from twilled Lie-Rinehart algebras\\
Banach Center Publications 51 (2000), 87--102
\endtitle
\author
Johannes Huebschmann
\endauthor
\date{November 29, 1998}
\enddate
\abstract {
Twilled L(ie-)R(inehart)-algebras generalize, in the Lie-Rinehart context, 
complex structures on smooth manifolds. An almost complex manifold determines 
an \lq\lq almost twilled pre-LR algebra\rq\rq, which is a true twilled 
LR-algebra iff the almost complex structure is integrable. We characterize 
twilled LR structures in terms of certain associated differential (bi)graded 
Lie and G(erstenhaber)-algebras; in particular the G-algebra arising from an 
almost complex structure is a d(ifferential) G-algebra iff the almost complex 
structure is integrable. Such G-algebras, endowed with a generator turning them
into a B(atalin-)V(ilkovisky)-algebra, occur on the B-side of the mirror 
conjecture. We generalize a result of Koszul to those dG-algebras which arise 
from twilled LR-algebras. A special case thereof explains the relationship 
between holomorphic volume forms and exact generators for the corresponding 
dG-algebra and thus yields in particular a conceptual proof of the 
Tian-Todorov
lemma. We give a differential homological algebra interpretation for twilled 
LR-algebras and by means of it we elucidate the notion of generator in terms of
homological duality for differential graded LR-algebras.}
\endabstract
\affil 
Universit\'e des Sciences et Technologies de Lille
\\
UFR de Math\'ematiques
\\
F-59 655 VILLENEUVE D'ASCQ C\'edex/France
\\
Johannes.Huebschmann\@univ-lille1.fr
\endaffil
\endkeywords
\endsubjclass
\endtopmatter
\document
\rightheadtext{Twilled Lie-Rinehart and BV algebras}
\medskip\noindent{\bf Introduction}\smallskip\noindent
In a series of seminal papers \cite\batviltw, \cite\batvilfo, \cite\batavilk, 
Batalin and Vilkovisky 
studied the quantization of constrained systems
and for that purpose introduced certain differential graded algebras
which have later been christened {\it Batalin-Vilkovisky algebras\/}.
Batalin-Vilkovisky algebras have recently become important
in string theory and elsewhere, cf. e.~g.
\cite\barakont,
\cite\getzltwo,
\cite\bv,
\cite\huebstas,
\cite\kosmathr,
\cite\liazutwo,
\cite\maninfiv,
\cite\stashnin,
\cite\xuone.
String theory leads to the
mysterious mirror conjecture.
A version thereof involves  
certain differential Batalin-Vilkovisky algebras
arising from a Calabi-Yau manifold.
A crucial ingredient is what is referred to
in the literature as the
{\it Tian-Todorov\/} lemma.
This Lemma, in turn, implies the

Here we will give a 
leisurely introduction
to a thorough study of
such differential Batalin-Vilkovisky
algebras and generalizations thereof
in the framework of Lie-Rinehart algebras,
thereby trying to avoid technicalities;
these and more details may be found in \cite\twilled.
\smallskip
A {\it Gerstenhaber\/} algebra is a graded commutative algebra
together with a bracket which (i) yields an ordinary
Lie bracket once the underlying module (or vector space)
has been regraded down by 1 and which (ii) satisfies
a certain derivation
property.
Such a bracket occurs in
Gerstenhaber's paper
\cite\gersthtw.
See Section 2 below for details.
A differential Batalin-Vilkovisky algebra 
is a Gerstenhaber algebra together with 
an exact generator, and the underlying Gerstenhaber algebras
of interest for us, in turn, arise 
as (bigraded) algebras of forms on {\it twilled Lie-Rinehart algebras\/}
(which we introduce below).
In the Lie-Rinehart context,
a twilled Lie-Rinehart algebra generalizes,
among others,
the notion of a complex structure
on a smooth manifold.
One of our results,
Theorem 2.3 below,
says that
 an \lq\lq almost
twilled Lie-Rinehart algebra\rq\rq\ 
is a true twilled Lie-Rinehart algebra
if and only if  
the corresponding
Gerstenhaber algebra is a differential Gerstenhaber algebra.
This implies, for example, that
the integrability condition for
an almost complex structure on a smooth manifold 
may be phrased as a condition saying that
a certain operator on
the corresponding Gerstenhaber algebra 
turns the latter into
a differential
Gerstenhaber algebra.
Now a theorem of Koszul 
\cite\koszulon\ 
establishes, on an ordinary smooth manifold,
a bijective correspondence between generators for 
the Gerstenhaber algebra of multi vector fields
and connections in the top exterior power of the tangent bundle
in such a way that
exact generators correspond to 
flat connections.
In Theorem 2.7 below
we  generalize this bijective correspondence
to the differential Gerstenhaber
algebras 
arising from twilled Lie-Rinehart algebras;
such Gerstenhaber algebras come into play,
for example, in the mirror conjecture.
What corresponds to a flat connection
on the line bundle in Koszul's theorem is now
a holomorphic volume form---its existence
is implied by the Calabi-Yau condition---and 
our generalization of Koszul's
theorem shows in particular how
a holomorphic volume form determines
a generator for the corresponding
differential Gerstenhaber algebra
turning it into a differential Batalin-Vilkovisky algebra.
The resulting 
differential Batalin-Vilkovisky algebra
then generalizes that which underlies what is called the B-model.
In particular,
as a consequence of our methods,
we obtain a new proof of the Tian-Todorov lemma.
We also give a 
differential homological algebra
interpretation of twilled Lie-Rinehart algebras
and, furthermore,
of a generator for a 
differential Batalin-Vilkovisky algebra
in terms of a suitable notion of homological duality.
This relies on results in our earlier papers
\cite\duality\ and \cite\bv\ as well as on various generalizations
therof.
\smallskip
I am indebted to Y. Kosmann-Schwarzbach and K. Mackenzie for
discussions, and to
J. Stasheff and A. Weinstein for some
e-mail correspondence
about various topics
related with the paper.
At the 
\lq\lq Poissonfest\rq\rq, 
Y. Kosmann-Schwarzbach
introduced me to the recent manuscript
\cite\schecone\ 
which treats topics somewhat related to
the present paper.
There is little overlap, though.
\smallskip
It is a pleasure to thank the organizers of the 
\lq\lq Poissonfest\rq\rq\
for the opportunity to present these results.
The material to be presented here is related to some of the work
of the late S. Zakrzewski; see Remark 4.3 below.
We respectfully  dedicate this 
paper to his memory.

\medskip
\noindent
{\bf 1. Twilled Lie-Rinehart algebras}
\smallskip\noindent
Let $R$ be a commutative ring.
A {\it Lie-Rinehart algebra\/} $(A,L)$ 
consists of a commutative $R$-algebra $A$ and an $R$-Lie algebra
$L$ together with 
an $A$-module structure
$A \otimes_R L \to L$ 
on $L$, written
$a \otimes_R\alpha \mapsto a\alpha$,
and an action
$L \to \roman{Der}(A)$ 
of $L$ on $A$ (which is a morphism of $R$-Lie algebras and)
whose adjoint 
$L \otimes_R A \to A$
is written
$\alpha\otimes_Ra \mapsto \alpha(a)$;
here
$a \in A$ and $\alpha \in L$.
These mutual actions 
are required
to satisfy certain compatibility properties modeled on
$(A,L) = (C^{\infty}(M),\roman{Vect}(M))$
where
$C^{\infty}(M)$ and $\roman{Vect}(M)$
refer to the algebra of smooth functions 
and the Lie algebra of smooth vector fields,
respectively,
on a smooth manifold $M$.
In general, the compatibility conditions read:
$$
\align
(a \alpha) b &= a \alpha (b),
\quad a, b \in A,\ \alpha \in L,
\tag1.1
\\
[\alpha, a \beta] &= \alpha (a) \beta+ a [\alpha, \beta],
\quad a \in A,\ \alpha,\beta \in L.
\tag1.2
\endalign
$$
For a Lie-Rinehart algebra $(A,L)$,
following \cite\rinehone,
we will refer to $L$ as an $(R,A)$-{\it Lie algebra\/}.
In differential geometry, $(R,A)$-Lie algebras arise as spaces of sections
of Lie algebroids.
Lie-Rinehart algebras have been studied
before Rinehart by
Herz \cite\herzone\ 
under the name \lq\lq pseudo-alg\`ebre de Lie\rq\rq\ 
as well as by Palais \cite\palaione\ 
under the name \lq\lq $d$-Lie-ring\rq\rq\.
We have chosen to refer to these object as Lie-Rinehart algebras
since Rinehart subsumed their
cohomology under standard homological algebra
and established a Poincar\'e-Birkhoff-Witt theorem for them \cite\rinehone.
In particular, Rinehart has shown how to describe de Rham cohomology
in the language of Ext-groups.
In a sense, 
the homological algebra interpretations of differential
Batalin-Vilkovisky algebras to be given below push
these observations of Rinehart's further.
\smallskip
Given two Lie-Rinehart algebras $(A,L')$ and $(A,L'')$, together with mutual
actions
$\cdot\colon L' \otimes_R L'' \to L''$
and
$\cdot\colon L'' \otimes_R L' \to L'$
which endow $L''$ and $L'$ with an
$(A,L')$- and
$(A,L'')$-module structure, respectively,
we will refer to
$(A,L',L'')$ as an {\it almost twilled Lie-Rinehart algebra\/};
we will call it a
{\it twilled Lie-Rinehart algebra\/}
provided the direct sum $A$-module structure on $L=L' \oplus L''$,
the sum $(L' \oplus L'')\otimes_R A \to A$
of the adjoints of the
$L'$- and $L''$-actions on $A$,
and the bracket $[\cdot,\cdot]$
on $L=L' \oplus L''$ given by
$$
[(\alpha'',\alpha'),(\beta'',\beta')]
=
[\alpha'',\beta''] +[\alpha',\beta']
+ \alpha'' \cdot \beta'
-\beta' \cdot \alpha''
+ \alpha' \cdot \beta''
-\beta'' \cdot \alpha'
\tag 1.3
$$
turn $(A,L)$
into a Lie-Rinehart algebra.
We then write
$L=L' \bowtie L''$ and refer to
$(A,L)$ as the
{\it twilled sum\/}
of
$(A,L')$ and
$(A,L'')$.
\smallskip
For illustration, consider
a smooth manifold $M$ with an almost complex structure,
let $A$ be the algebra of smooth complex functions on $M$,
$L$ the $(\Bobb C,A)$-Lie algebra of complexified smooth vector fields
on $M$, and consider the ordinary decomposition
of the complexified tangent bundle $\tau^{\Bobb C}_M$ 
as a direct sum 
$\tau'_M \oplus \tau''_M$ 
of the {\it almost holomorphic\/}
and {\it almost antiholomorphic\/}
tangent bundles
$\tau'_M$
and  $\tau''_M$, respectively;
write $L'$ and $L''$ for their spaces of smooth sections.
Then 
$(A,L',L'')$,
together with the mutual actions coming from $L$,
is a twilled Lie-Rinehart algebra
if and only if the almost complex structure is integrable,
i.~e. a true complex structure;
$\tau'_M$
and  $\tau''_M$
are then the ordinary {\it holomorphic\/}
and {\it antiholomorphic\/}
tangent bundles, respectively.
The precise analogue of an almost complex structure
is what we call
an {\it almost twilled pre-Lie-Rinehart algebra\/} structure;
this notion is weaker than that of 
almost twilled Lie-Rinehart algebra.
The basic difference is that,
for an almost twilled pre-Lie-Rinehart algebra,
instead of having mutual
actions
$\cdot\colon L' \otimes_R L'' \to L''$
and
$\cdot\colon L'' \otimes_R L' \to L'$,
we only require that there be given
$R$-linear pairings
$\cdot\colon L' \otimes_R L'' \to L''$
and
$\cdot\colon L'' \otimes_R L' \to L'$,
which endow $L''$ and $L'$ with an
$(A,L')$- and
$(A,L'')$-connection, respectively;
see \cite \twilled\  for details.
A situation similar to that of a complex
structure on a smooth manifold
and giving rise to a twilled Lie-Rinehart algebra
arises from a smooth manifold with two
transverse foliations as well as from a Cauchy-Riemann structure
(cf. \cite\canhawei);
see \cite\twilled\ for some comments about Cauchy-Riemann structures.
Lie bialgebras
provide
another class of examples
of
twilled Lie-Rinehart algebras; 
Kosmann-Schwarzbach and F. Magri refer to these objects, or rather
to the corresponding twilled sum, as
{\it twilled extensions of Lie algebras\/} \cite\kosmagtw;
Lu and Weinstein 
call them
{\it double Lie algebras\/} \cite\luweinst;
and Majid
uses the terminology
{\it matched pairs\/}
of Lie algebras \cite\majidtwo.
Spaces of sections
of suitable pairs of Lie algebroids
with additional structure
lead to yet another class of 
examples of
twilled Lie-Rinehart algebras;
these have been studied in the literature
under the name
{\it matched pairs 
of Lie algebroids\/}
by Mackenzie \cite\mackfift\  and Mokri \cite\mokritwo.
\smallskip
An almost twilled Lie-Rinehart algebra $(A,L'',L')$
is a true 
twilled Lie-Rinehart algebra 
if and only if 
$(A,L'',L')$
satisfies three compatibility conditions;
these are spelled out in \cite\twilled\ (Proposition 1.7).
This proposition is merely an adaption of
earlier results in the literature
to our more general situation.
Another interpretation of the compatibility conditions
involves
certain annihilation properties of the two operators
$d'$ and $d''$
which arise,
given
an almost twilled pre-Lie-Rinehart algebra
$(A,L',L'')$,
as formal extensions of the ordinary Lie-Rinehart
differentials
with respect to $L'$ and $L''$, respectively,
on the bigraded algebra
$\roman{Alt}_A^*(L'', \roman{Alt}_A^*(L',A))$
(but are not necessarily exact):
To explain this, we will say
that an $A$-module $M$ has {\it property\/} P 
provided for  $x\in M$,
$\phi(x) = 0$ for every $\phi \colon M \to A$
implies that $x$ is zero.
For example, a projective $A$-module has property P.
We now have the following,
cf. Theorem 1.15 in \cite\twilled.

\proclaim{Theorem 1.4}
If $(A,L',L'')$ is a twilled Lie-Rinehart algebra,
the operators $d'$ and $d''$ turn
the bigraded algebra
$\roman{Alt}_A^*(L'', \roman{Alt}_A^*(L',A))$
into a differential bigraded algebra
which
then necessarily computes
the cohomology
$\roman H^*(\roman{Alt}_A(L,A))$
of the twilled sum $L$ of $L'$ and $L''$.
Conversely, given
an almost twilled pre-Lie-Rinehart algebra
$(A,L',L'')$,
if 
the operators $d'$ and $d''$ turn
the bigraded algebra
$\roman{Alt}_A^*(L'', \roman{Alt}_A^*(L',A))$
into a differential bigraded algebra
and if $L'$ and $L''$ have property {\rm P},
$(A,L',L'')$ is a true twilled Lie-Rinehart algebra.
\endproclaim

For example,
{\sl for the twilled Lie-Rinehart algebra
arising from the holomorphic and antiholomorphic
tangent bundles of a complex manifold,
the resulting differential
bigraded algebra\/}
$(\roman{Alt}_A^*(L'', \roman{Alt}_A^*(L',A)),\parti', \parti'')$
{\sl comes down to the ordinary Dolbeault complex\/}.

\medskip
\noindent
{\bf 2. Twilled Lie-Rinehart algebras, Gerstenhaber algebras, and dBV-algebras}
\smallskip
\noindent
We now explain briefly how
other characterizations
of
twilled Lie-Rinehart algebras
explain 
the structure of certain differential Batalin-Vilkovisky
algebras.
Section 2 of \cite\twilled\  is devoted to 
more details about
differential graded Lie-Rinehart
algebras.
\smallskip
Let $(A,L)$ be a Lie-Rinehart algebra, and let
$\Cal A$ be a graded commutative 
$A$-algebra 
which is endowed with a graded $(A,L)$-module structure
in such a way that (i) $L$ acts
on $\Cal A$ by derivations---this is equivalent
to requiring the structure map
from $\Cal A \otimes_A \Cal A$ to $\Cal A$
to be a morphism of graded $(A,L)$-modules---and that (ii)
the canonical map from $A$ to $\Cal A$ is
a morphism of left $(A,L)$-modules.
Let 
$\Cal L = \Cal A \otimes_A L$, and define a bigraded bracket
$$
[\cdot,\cdot]
\colon
\Cal L
\otimes_R
\Cal L
@>>>
\Cal L
\tag2.1.1
$$
of bidegree $(0,-1)$
by means of the formula
$$
[\alpha \otimes_A x, \beta \otimes_A y]
=
(\alpha \beta) \otimes_A[x,y]
+\alpha (x\cdot \beta) \otimes_A y
- (-1)^{|\alpha||\beta|} \beta (y \cdot \alpha)\otimes_A x
\tag2.1.2
$$
where
$\alpha,\beta \in \Cal A$ and $x,y \in L$.
A calculation shows that,
for every $\beta \in \Cal A$ and every $x,y,z \in L$,
$$
[[x,y],\beta \otimes_Az] -
([x,[y,\beta \otimes_Az]] -[y,[x,\beta \otimes_A z]])
=\big([x,y] (\beta) - x(y(\beta))-y(x(\beta))\big )\otimes_A z,
$$
whence 
(2.1.1) being a graded Lie bracket is actually equivalent
to the structure map $L\otimes_R\Cal A \to \Cal A$ being a Lie algebra
action.
Moreover,
let
$$
\Cal A \otimes_R \Cal L @>>> \Cal L
\tag2.1.3
$$
be the obvious graded left
$\Cal A$-module structure
arising from extension of scalars,
that is from extending $L$ to a (graded)
$\Cal A$-module,
and
define a pairing
$$
\Cal L \otimes_R \Cal A @>>> \Cal A
\tag2.1.4
$$
by
$$
(\alpha \otimes_A x) \otimes_R \beta \mapsto
(\alpha \otimes_A x) (\beta)
=\alpha (x(\beta)).
\tag2.1.5
$$
Then $(\Cal A,\Cal L)$,
together with (2.1.1),
(2.1.3) and (2.1.4),
constitutes a graded Lie-Rinehart algebra.
We refer to 
$(\Cal A,\Cal L)$
as the {\it crossed product\/}
of  $\Cal A$ and $(A,L)$
and to the corresponding
$(R,\Cal A)$-Lie algebra $\Cal L$
as the {\it crossed product\/} of $\Cal A$ and $L$.
More details about this notion of graded crossed product
Lie-Rinehart algebra may be found in \cite\twilled\ (2.8).
\smallskip
\noindent
{\smc Remark 2.1.6.} We must be a little circumspect here:
The three terms on the right-hand side of (2.1.2)
are {\it not\/} well defined individually; only their sum is well
defined. For example,
if we take $ax$ instead of $x$, where $a \in A$,
on the left-hand side,
$\alpha \otimes_A(ax)$ equals
$(\alpha a)\otimes_A x$
but
$(\alpha \beta) \otimes_A[ax,y]$
differs from
$(\alpha a\beta) \otimes_A[x,y]$.
\smallskip
Let
$(A,L'',L')$ be an almost twilled Lie-Rinehart algebra
having $L'$ finitely generated and projective as an $A$-module.
Write $\Cal A'' = \roman{Alt}_A(L'',A)$
and
$\Cal L' = \roman{Alt}_A(L'',L')$.
Now
$\Cal A''$ is a graded commutative $A$-algebra
and, endowed with the Lie-Rinehart differential $d''$
(which corresponds to the
$(R,A)$-Lie algebra structure on $L''$),
$\Cal A''$ is a differential graded commutative $R$-algebra.
Moreover,
from the $(A,L'')$-module structure
on $L'$,
$\Cal L'$
inherits an obvious
differential graded
$\Cal A''$-module structure.
Furthermore,
the $(A,L')$-structure
on $L''$ induces an action
of $L'$ on
$\Cal A''$
by graded derivations.
Since $L'$ is supposed to be finitely generated and projective
as an $A$-module,
the canonical $A$-module morphism
$$
\Cal A''
\otimes_A
L
@>>>
\Cal L'
=
\roman{Alt}_A(L'',L')
$$
is an isomorphism of graded $A$-modules,
in fact of graded
$\Cal A''$-modules.
Taking $L= L'$ and $\Cal A = \Cal A''$,
the graded crossed product Lie-Rinehart
structure explained above is available,
and we have
the graded crossed product
Lie-Rinehart algebra
$(\Cal A'',\Cal L')$.
Now the
$(R,A)$-Lie algebra structure on $L''$ and
the $(A,L'')$-module structure on $L'$
determine the corresponding Lie-Rinehart differential
on $\Cal L' = \roman{Alt}_A(L'',L')$;
we denote it by $d''$.
By symmetry, 
when $L''$ is finitely generated and projective as an $A$-module,
we have the same structure, with
$L'$ and $L''$ interchanged.

\proclaim{Theorem 2.1}
As an $A$-module, $L'$ being supposed to be finitely generated and projective,
the statements
{\rm (i)}, {\rm (ii)}, and {\rm (iii)} below are equivalent:
\newline\noindent
{\rm (i)}
$(A,L'',L')$
is a true
twilled Lie-Rinehart algebra;
\newline\noindent
{\rm (ii)}
$(\Cal L',d'') = (\roman{Alt}_A(L'',L'),d'')$
is a differential graded $R$-Lie algebra;
\newline\noindent
{\rm (iii)}
$(\Cal A'',\Cal L';d'')$ is a differential graded Lie-Rinehart
algebra.
\newline\noindent
Thus, under these circumstances,
there is a bijective correspondence between
twilled Lie-Rinehart algebra and differential graded
Lie-Rinehart algebra structures.
\endproclaim
For a proof of this result
and for more details, see  (3.2) in \cite\twilled.
We note that, in the situation of Theorem 2.1,
the Lie bracket on
$\Cal L'=\roman{Alt}_A(L'',L')$
does not just come down to the shuffle product of
forms on $L''$ and the Lie bracket on $L'$;
in fact, such a bracket would not even be well defined
since the Lie bracket of 
$L'$ 
is not $A$-linear, i.~e.,
in the usual differential geometry context,
does not behave as a \lq\lq tensor\rq\rq.
\smallskip
In \cite\crosspro, a notion of graded crossed product
Lie-Rinehart algebra has been given
which may entirely be described in terms of
the $A$-modules of forms $\roman{Alt}_A(L'',A)$
and $\roman{Alt}_A(L'',L')$;
by means thereof, the statement of
Theorem 2.1
may be established without the hypothesis that
$L'$ be finitely generated and projective as an $A$-module.
\smallskip
When $(A,L',L'')$
is the twilled Lie-Rinehart algebra
arising from 
the holomorphic and antiholomorphic tangent bundles of
a smooth complex manifold $M$,
$(\Cal L',d'') = (\roman{Alt}_A(L'',L'),d'')$
is what is called the
{\it Kodaira-Spencer\/}
algebra in the literature;
it controls the infinitesimal deformations
of the complex structure on $M$.
The cohomology
$\roman H^*(L'',L')$
then inherits a graded Lie algebra structure
and the obstruction to deforming the complex structure
is the map
$\roman H^1(L'',L')\to \roman H^2(L'',L')$
which sends $\eta\in \roman H^1(L'',L')$ to 
$[\eta,\eta]\in \roman H^2(L'',L')$.
\smallskip
Recall that a {\it Gerstenhaber algebra\/}
is a graded commutative $R$-algebra
$\Cal A$
together with a graded Lie bracket
from $\Cal A \otimes_R \Cal A$ to $\Cal A$
of degree $-1$ (in the sense that,
if $\Cal A$ is regraded down by one,
$[\cdot,\cdot]$ is an ordinary graded Lie bracket)
such that, for each homogeneous element $a$ of $\Cal A$,
$[a,\cdot]$ is a derivation of $\Cal A$ of degree $|a|-1$
where $|a|$ refers to the degree of $a$;
see \cite\geschthr\ 
where these objects are called G-algebras,
or \cite{\bv,\kosmathr,\liazutwo,\xuone}.
In our paper \cite\bv,
we worked out an intimate link between
Gerstenhaber's and Rinehart's papers \cite\gersthtw\ 
and \cite\rinehone\ 
which involves the notion of Gerstenhaber bracket.
In a sense,
in the present paper we extend this link to the differential
graded situation.
\smallskip
Given  a bigraded commutative $R$-algebra
$\Cal A$,
we will say that
a bigraded bracket
$[\cdot,\cdot]\colon\Cal A \otimes_R \Cal A \to \Cal A$
of bidegree $(0,-1)$
is a {\it bigraded\/}
{\it Gerstenhaber bracket\/}
provided 
$[\cdot,\cdot]$ is an ordinary bigraded Lie bracket
when the second degree of
$\Cal A$ is regraded down by one, the first one being kept,
such that, for each homogeneous element $a$ of $\Cal A$
of bidegree $(p,q)$,
$[a,\cdot]$ is a derivation of $\Cal A$ of bidegree $(p,q-1)$;
a bigraded
$R$-algebra
with a
bigraded Gerstenhaber
bracket will be referred to
as
a {\it bigraded Gerstenhaber algebra\/}.
Moreover, given
a bigraded Gerstenhaber algebra
$(\Cal A,[\cdot,\cdot])$
together with a 
differential $d$
of bidegree $(1,0)$
which endows $\Cal A$ with
a differential graded
$R$-algebra structure
we will say that
$(\Cal A,[\cdot,\cdot])$
and $d$ constitute
a {\it differential bigraded\/} Gerstenhaber algebra
(or differential bigraded G-algebra),
written $(\Cal A,[\cdot,\cdot], d)$, 
provided
$d$ behaves as a derivation for
the bigraded Gerstenhaber bracket $[\cdot,\cdot]$, that is,
$$
d[x,y] = [dx,y] -(-1)^{|x|} [x,dy], \quad x,y \in \Cal A,
$$
where the total degree $|x|$ is the sum of the bidegrees.
\smallskip
Recall that, given a Lie-Rinehart algebra $(A,L)$,
the Lie bracket on $L$ determines 
a Gerstenhaber bracket
on the exterior $A$-algebra
$\Lambda_AL$ on $L$; 
for $\alpha_1,\ldots,\alpha_n \in L$,
the bracket $[u,v]$ in $\Lambda_AL$
of 
$u=\alpha_1\wedge \ldots \wedge\alpha_{\ell}$ and
$v=\alpha_{\ell+1}\wedge \ldots \wedge \widehat{\alpha_n}$
is given by the expression
$$
[u,v]
=
(-1)^{\ell} 
\sum_{j\leq \ell <k} (-1)^{(j+k)}
\lbrack \alpha_j,\alpha_k \rbrack \wedge 
\alpha_1\wedge \ldots \widehat{\alpha_j} \ldots \widehat{\alpha_k}
\ldots \wedge \alpha_n,
\tag2.2.1
$$
where $\ell = |u|$ is the degree of $u$,
cf. \cite\bv\ (1.1).
\smallskip
We now return  to a general almost twilled Lie-Rinehart algebra
$(A,L',L'')$ having $L'$ finitely generated and projective
as an $A$-module and
consider the graded crossed product
Lie-Rinehart algebra
$(\Cal A'',\Cal L')$. 
The graded Lie-Rinehart bracket on $\Cal L' (=\roman{Alt}_A(L'',L'))$
extends to a (bigraded) bracket
on
$\roman{Alt}_A(L'',\Lambda_AL')$
which turns the latter into a
bigraded Gerstenhaber algebra;
as a bigraded algebra,
$\roman{Alt}_A(L'',\Lambda_AL')$
could be thought as of the exterior
$\Cal A''$-algebra on
$\Cal L'$, and we write sometimes
$$
\Lambda_{\Cal A''} \Cal L' = \roman{Alt}_A(L'',\Lambda_AL').
$$
With reference to the
graded Lie bracket
$[\cdot,\cdot]$
on $\Cal L'$
and the $L'$-action on $\Cal A''$,
the bigraded Gerstenhaber bracket
$$
[\cdot,\cdot]
\colon
\Lambda_{\Cal A''}\Cal L'
\otimes_R
\Lambda_{\Cal A''}\Cal L'
@>>>
\Lambda_{\Cal A''}\Cal L'
\tag2.2.2
$$
on
$\Lambda_{\Cal A''}\Cal L'$
may be described
by the formulas
$$
\aligned
[\alpha \beta,\gamma] 
&=
\alpha [\beta,\gamma] +{(-1)}^{|\alpha||\beta|} \beta[\alpha,\gamma],
\quad \alpha,\beta,\gamma \in \Lambda_{\Cal A''}\Cal L',
\\
[x,a] &= x(a), \quad x \in L',\, a \in\Cal A'',
\\ 
[\alpha,\beta] 
&= -{(-1)}^{(|\alpha|-1)(|\beta|-1)}[\beta,\alpha],\quad
\alpha,\beta \in \Lambda_{\Cal A''}\Cal L',
\endaligned
\tag2.2.3
$$
where $|\cdot|$ refers to the total degree, i.e. the sum of the
two bidegree components.
The bracket (2.2.2)
is in fact the {\it (bigraded)
crossed product bracket extension\/}
of the Gerstenhaber bracket
on $\Lambda_AL'$
and $\Lambda_{\Cal A''}\Cal L'$
may be viewed as the {\it (bigraded) crossed product Gerstenhaber algebra\/}
of $\Cal A''$ with
the ordinary Gerstenhaber algebra
$\Lambda_A L'$.
See Section 4 of \cite\twilled\ for details.
\smallskip
The Lie-Rinehart differential $d''$
which corresponds to the
Lie-Rinehart structure on $L''$
and the induced
graded $(A,L'')$-module structure
on $\Lambda_AL'$
turn
$\roman{Alt}_A(L'',\Lambda_AL')$
into a differential (bi)-graded commutative $R$-algebra.
By symmetry, 
when $L''$ is finitely generated and projective as an $A$-module,
we have the same structure, with
$L'$ and $L''$ interchanged.

\proclaim{Theorem 2.3}
The almost twilled Lie-Rinehart algebra $(A,L'',L')$
is a true
twilled Lie-Rinehart algebra
if and only if
$(\Lambda_{\Cal A''}\Cal L',d'')$
($=(\roman{Alt}_A(L'',\Lambda_AL'),d'')$)
is a differential (bi)-graded Gerstenhaber algebra.
\endproclaim

See Theorem 4.4 in \cite\twilled\ and its proof.
\smallskip
When $(A,L',L'')$ 
arises from 
the holomorphic and antiholomorphic tangent bundles of
a smooth complex manifold $M$,
the resulting differential Gerstenhaber algebra
$(\roman{Alt}_A(L'',\Lambda_AL'),d'')$
is that of forms of type
$(0,*)$
with values in the holomorphic multi vector fields,
the operator $d''$
being the Cauchy-Riemann operator
(which is more usually written $\overline \partial$).
This 
differential Gerstenhaber algebra
comes into play in the mirror conjecture;
it was studied by
Barannikov-Kontsevich \cite\barakont,
Manin \cite\maninfiv,
Witten \cite\wittetwe,
and others.
\smallskip
Let now $(A,L'',L')$ be a twilled Lie-Rinehart algebra
having $L'$
finitely generated and projective as an $A$-module
of constant rank $n$ (say),
and write
$\Lambda_A^nL'$ for the top exterior power of $L'$ over $A$.
Consider the differential Gerstenhaber algebra
$(\roman{Alt}_A(L'',\Lambda_AL'),d'')$.
Our next aim is to study
generators thereof.
To this end, we observe that,
when
$\roman{Alt}_A(L',\Lambda_A^nL')$
is endowed with the obvious graded
$(A,L'')$-module structure
induced from the left $(A,L'')$-module structure
on $L'$ which is part of the structure of twilled Lie-Rinehart algebra,
the canonical isomorphism
$$
\roman{Alt}_A(L'',\Lambda_AL')
@>>>
\roman{Alt}_A(L'',\roman{Alt}_A(L',\Lambda_A^nL'))
\tag2.4
$$
of graded $A$-modules
is compatible with the differentials
which correspond to the Lie-Rinehart structure
on $L''$ and the
$(A,L'')$-module structures on the coefficients
on both sides of (2.4);
abusing notation, we denote each of these differentials by $d''$.

\smallskip
For a bigraded Gerstenhaber algebra $\Cal A$
over $R$,
with bracket operation written $[\cdot,\cdot]$,
an $R$-linear operator $\Delta$
on $\Cal A$  
of bidegree $(0,-1)$
will be said to {\it generate\/}
the Gerstenhaber bracket
provided, for every homogeneous $a, b \in \Cal A$,
$$
[a,b] = (-1)^{|a|}\left(
\Delta(ab) -(\Delta a) b - (-1)^{|a|} a (\Delta b)\right);
\tag2.5
$$
the operator $\Delta$ is then called a {\it generator\/}.
A generator $\Delta$  is said to be {\it exact\/}
provided $\Delta\Delta$ is zero, that is, $\Delta$ is a differential;
an exact generator will henceforth be written $\partial$.
A bigraded Gerstenhaber algebra $\Cal A$
together with a generator $\Delta$
will be 
called a {\it weak bigraded Batalin-Vilkovisky\/} algebra
(or weak bigraded BV-algebra);
when
the generator is exact,
written $\partial$,
we will  refer to 
$(\Cal A,\partial)$
(more simply)
as a
{\it bigraded Batalin-Vilkovisky\/} algebra
(or bigraded BV-algebra).
\smallskip
It is clear that a generator
determines the bigraded Gerstenhaber bracket.
An observation due to Koszul \cite\koszulon\ 
(p. 261)
carries over to the bigraded case:
for any  bigraded Batalin-Vilkovisky algebra
$(\Cal A, [\cdot,\cdot], \partial)$,
the operator $\partial$ 
(which is exact by assumption)
behaves as a derivation for
the bigraded Gerstenhaber bracket $[\cdot,\cdot]$,
that is,
$$
\partial [x,y] = [\partial x,y] -(-1)^{|x|} [x,\partial y], 
\quad x,y \in \Cal A.
\tag2.6
$$
An exact generator $\partial$
does in general {\it not\/} behave
as a derivation for the multiplication
of $\Cal A$, though.
Let $(\Cal A,\Delta)$ be a weak bigraded Batalin-Vilkovisky algebra, 
write $[\cdot,\cdot]$
for the bigraded Gerstenhaber bracket generated by $\Delta$,
and let $d$ be a differential of bidegree $(+1,0)$
which endows $(\Cal A,[\cdot,\cdot])$ with a differential bigraded 
Gerstenhaber algebra structure. Consider the graded commutator
$[d,\Delta] = d\Delta + \Delta d$ 
on $\Cal A$;
it is an 
operator of bidegree $(1,-1)$ and hence of total degree zero.
We will say that
$(\Cal A,\Delta,d)$
is a {\it weak differential\/}
bigraded Batalin-Vilkovisky algebra provided the commutator
$[d,\Delta]$
is zero.
In particular,
a weak differential
bigraded Batalin-Vilkovisky algebra 
$(\Cal A,\partial,d)$
which has $\partial$ exact
will be called a 
{\it differential bigraded Batalin-Vilkovisky algebra\/}.
On 
the underlying bigraded object $\Cal A$ of a differential
bigraded Batalin-Vilkovisky algebra 
$(\Cal A,\partial,d)$,
the graded commutator
$[d,\partial]$ 
manifestly behaves as a derivation for the 
bigraded
Gerstenhaber bracket
since
$d$ and $\partial$
both
behave as derivations for this bracket.

We now reproduce the statement of Theorem 5.4.6 in \cite\twilled.

\proclaim{Theorem 2.7}
The isomorphism {\rm (2.4)} furnishes a bijective correspondence
between generators
of the bigraded Gerstenhaber structure
on the left-hand side of {\rm (2.4)}  and
$(A,L')$-connections
on $\Lambda_A^nL'$
in such a way that
exact generators
correspond to
$(A,L')$-module structures
(i.~e. flat connections).
Under this correspondence,
generators
of the differential bigraded Gerstenhaber structure
on the left-hand side 
correspond to
$(A,L')$-connections
on $\Lambda_A^nL'$
which are compatible with the
$(A,L'')$-module structure on $\Lambda_A^nL'$.
\endproclaim

Thus, in particular,
exact generators
of the differential bigraded Gerstenhaber structure
on the left-hand side 
correspond to
$(A,L'')$-compatible $(A,L')$-module structures
on $\Lambda_A^nL'$.
\smallskip
When $L''$ is trivial and
$L'$ the Lie algebra of smooth vector fields
on a smooth manifold,
the statement of this theorem comes down to the result
of Koszul \cite\koszulon\ mentioned earlier.
Our result not only provides many examples 
of differential Batalin-Vilkovisky algebras
but also explains how
every differential Batalin-Vilkovisky algebra 
having an underlying bigraded $A$-algebra of the kind
$\roman{Alt}_A(L'',\Lambda_AL')$
arises.
\smallskip
When $(A,L',L'')$ is the twilled Lie-Rinehart algebra
which comes from the holomorphic and antiholomorphic tangent bundles
of a smooth complex manifold $M$ as explained above,
the theorem gives a bijective correspondence
between 
generators
of the 
differential bigraded Gerstenhaber algebra
$(\roman{Alt}_A(L'',\Lambda_AL'),d'')$
of forms of type
$(0,*)$
with values in the holomorphic multi vector fields,
the differential $d''$ being the
Cauchy-Riemann operator,
and holomorphic connections
on the highest exterior power of the holomorphic tangent bundle
in such a way that
exact generators correspond to flat holomorphic connections.
In particular, 
suppose
that $M$ is a {\it Calabi-Yau\/}
manifold,
that is,
admits
a holomorphic volume form $\Omega$ (say).
This holomorphic volume form
identifies
the highest exterior power of the holomorphic tangent bundle
with the algebra of 
smooth complex functions on $M$
as a module over $L=L''\oplus L'$,
hence induces
a flat holomorphic connection
thereupon 
and thence
an exact generator
$\partial_\Omega$
for
$(\roman{Alt}_A(L'',\Lambda_AL'),d'')$,
turning the latter into a differential (bi)graded
Batalin-Vilkovisky algebra.
This is precisely the
differential (bi)graded
Batalin-Vilkovisky algebra
coming into play 
on the B-side of the mirror conjecture
and studied in the cited sources.
The fact that the holomorphic volume form induces a generator
for the differential Gerstenhaber structure
is referred to in the literature as
the {\it Tian-Todorov lemma\/};
cf \cite\tianone,
\cite\todortwo.
In our approach, 
this lemma drops out
as a special case of our generalization
of Koszul's theorem to the bigraded setting,
and this generalization indeed provides a conceptual proof of the lemma.
This lemma implies that, for a K\"ahlerian Calabi-Yau manifold $M$,
the deformations of the complex structure are unobstructed, that is to say,
there is an open subset of
$\roman H^1(M,\tau_M)$
parametrizing the deformations of the complex structure;
here  $\roman H^1(M,\tau_M)$
is the first cohomology group
of $M$ with values in the holomorphic tangent bundle
$\tau_M$.
Under these circumstances,
after a choice of holomorphic volume form
$\Omega$
has been made,
the canonical isomorphism
(2.4), combined
with the isomorphism
$$
\Omega^{\flat}\colon
\roman{Alt}_A(L'',\roman{Alt}_A(L',\Lambda_A^nL'))
@>>>
\roman{Alt}_A(L'',\roman{Alt}_A(L',A))
$$
induced by $\Omega$
identifies
$(\roman{Alt}_A(L'',\Lambda_AL'),d'',\partial_\Omega)$
with the Dolbeault complex
of $M$ and hence
the cohomology
$\roman H^*(\roman{Alt}_A(L'',\Lambda_AL'),d'',\partial_\Omega)$
with the ordinary complex valued cohomology of $M$.
This is nowadays well understood.
The cohomology
$\roman H^*(\roman{Alt}_A(L'',\Lambda_AL'),d'',\partial_{\Omega})$
is referred to 
in the literature as the {\it extended moduli space of complex
structures\/} \cite\wittetwe;
it underlies what is called the B-{\it model\/}
in the theory of mirror symmetry.

\medskip
\noindent
{\bf 3. Twilled Lie-Rinehart algebras and differential homological algebra}
\smallskip\noindent
We now spell out interpretations of some of the above
results in terms of differential homological algebra.
\smallskip
Let
$(A,L',L'')$ be  a twilled Lie-Rinehart algebra
having $L'$ and $L''$
finitely generated and projective as $A$-modules.
Let $(\Cal A'',\Cal L';d'')$ 
be the differential graded
crossed product Lie-Rinehart algebra
$(\roman{Alt}_A(L'',A),\roman{Alt}_A(L'',L');d'')$
mentioned before.
Let
$L=L' \bowtie L''$ be the twilled sum
of $L'$ and $L''$, and
consider 
the differential graded Lie-Rinehart cohomology
$\roman H^*(\Cal L',\Cal A'')$.
see Section 6 in \cite\twilled\ for details
where also a proof of the following result may be found;
cf. in particular Theorem 6.15 in \cite\twilled.

\proclaim{Theorem 3.1}
The differential
bigraded algebra
$(\roman{Alt}_A^*(L'', \roman{Alt}_A^*(L',A)),\parti', \parti'')$
computes
the differential graded Lie-Rinehart cohomology
$\roman H^*(\Cal L',\Cal A'')$.
Moreover,
this
differential graded Lie-Rinehart cohomology
is naturally isomorphic to
the Lie-Rinehart cohomology 
$\roman H^*(L,A)$.
\endproclaim

When $L''$ is trivial,
so that 
$\roman H^*(\Cal L',\Cal A'')$
is an ordinary (ungraded) Lie-Rinehart algebra
$(A,L)$,
the differential graded Lie-Rinehart cohomology
boils down to 
the ordinary  Lie-Rinehart cohomology
$\roman H^*(L,A)$.
Moreover, 
for the special case 
when
$A$ and $L$ are the algebra of smooth functions and smooth vector
fields on a smooth manifold,
the Lie-Rinehart cohomology
$\roman H^*(L,A)$
amounts to the de Rham cohomology;
this fact has been established by Rinehart \cite\rinehone.
In our more general situation,
when
the twilled Lie-Rinehart algebra
$(A,L',L'')$ 
arises from the holomorphic and antiholomorphic tangent bundles
of a smooth complex manifold,
the complex
calculating
the differential graded Lie-Rinehart cohomology
$\roman H^*(\Cal L',\Cal A'')$
of the differential graded
crossed product Lie-Rinehart algebra
$(\Cal A'',\Cal L';d'')= (\roman{Alt}_A(L'',A),\roman{Alt}_A(L'',L');d'')$
is the Dolbeault complex,
and 
the differential graded Lie-Rinehart cohomology
amounts to the Dolbeault cohomology.
Thus our approach provides, in particular, 
an interpretation of the Dolbeault complex in the framework
of differential homological algebra.
\smallskip
Generalizing results in our earlier paper \cite\bv, 
we can now elucidate
the concept of generator of a differential bigraded
Batalin-Vilkovisky algebra
in the framework of homological duality for differential graded
Lie-Rinehart algebras
in the following way:
{\sl An exact generator amounts to the differential
in a standard complex computing differential graded Lie-Rinehart
homology\/} (!) 
{\sl with appropriate coefficients\/};
see Proposition 7.13 in \cite\twilled\   for details.
It may then be shown that,
when the appropriate additional structure 
(in terms of Lie-Rinehart differentials and dBV-generators)
is taken into account,
the above isomorphism (2.4) is essentially just a duality
isomorphism 
in the (co)homology of the differential graded crossed
product Lie-Rinehart algebra
$(\Cal A'',\Cal L')$;
see  Proposition 7.14 in \cite\twilled\ for details.
In particular,
the Tian-Todorov Lemma comes down to a statement
about differential graded
(co)homological duality.
\medskip\noindent
{\bf 4. Twilled Lie-Rinehart algebras and Lie-Rinehart bialgebras}
\smallskip\noindent
Twilled Lie-Rinehart algebras thus generalize Lie bialgebras,
and the twilled sum
is an analogue,
even a generalization, of the Manin double of a Lie bialgebra.
The 
Lie bialgebroids
introduced by Mackenzie and Xu 
\cite\mackxu\ generalize Lie bialgebras as well,
and there is a corresponding notion of Lie-Rinehart bialgebra.
However,
twilled Lie-Rinehart algebras and Lie-Rinehart bialgebras
are different, in fact non-equivalent notions
which both generalize
Lie bialgebras.
In a sense,
Lie-Rinehart bialgebras
generalize Poisson 
and in particular symplectic
structures
while twilled Lie-Rinehart algebras
generalize complex structures.
We now give a
a characterization of twilled Lie-Rinehart algebras
in terms of Lie-Rinehart bialgebras.
See Theorem 4.8 in \cite\twilled\ 
for more details.
Let $L$ and $D$ be $(R,A)$-Lie algebras
which,
as $A$-modules, are  finitely generated and projective,
in such a way that,
as an $A$-module, $D$ is isomorphic to $L^* = \roman{Hom}_A(L,A)$.
We say
that $L$ and $D$ {\it are in duality\/}.
We write $d$ for the differential on
$\roman{Alt}_A(L,A)\cong \Lambda_AD$
coming from the Lie-Rinehart structure on $L$
and
$d_*$ for the differential on
$\roman{Alt}_A(D,A)\cong \Lambda_AL$
coming from the Lie-Rinehart structure on $D$.
Likewise we denote 
the Gerstenhaber bracket
on
$\Lambda_A L$
coming from 
the Lie-Rinehart structure on 
$L$ by
$[\cdot,\cdot]$
and
that
on
$\Lambda_A D$
coming from 
the Lie-Rinehart structure on 
$D$ by
$[\cdot,\cdot]_*$.
We will say that $(A,L,D)$
constitutes a
{\it Lie-Rinehart\/}
bialgebra if the differential $d_*$ on
$\roman{Alt}_A(D,A) \cong \Lambda_A L$
and the Gerstenhaber bracket
$[\cdot,\cdot]$
on
$\Lambda_A L$
are related by
$$
d_*[x,y] = [d_*x,y] + [x,d_*y], \quad x,y \in L,
$$
or equivalently,
if the differential $d$ on
$\roman{Alt}_A(L,A) \cong \Lambda_A D$
behaves as a derivation for
the Gerstenhaber bracket
$[\cdot,\cdot]_*$
in all degrees, that is to say
$$
d[x,y]_* = [dx,y]_* -(-1)^{|x|} [x,dy]_*, \quad x,y \in \Lambda_AD.
$$
Thus, for a Lie-Rinehart bialgebra $(A,L,D)$,
$$
(\Lambda_A L, [\cdot,\cdot], d_*)
=
(\roman{Alt}_A(D,A), [\cdot,\cdot], d_*)
$$
is a strict {\it differential Gerstenhaber algebra\/},
and the same is true of
$$
(\Lambda_A D, [\cdot,\cdot]_*, d)
=
(\roman{Alt}_A(L,A), [\cdot,\cdot]_*, d);
$$
see \cite\kosmathr\ (3.5) for details.
In fact, a straightforward extension of an observation of
Y. Kosmann-Schwarzbach \cite\kosmathr\ 
shows that
Lie-Rinehart bialgebra structures
on $(A,L,D)$
and
strict differential Gerstenhaber algebra structures
on
$(\Lambda_A L, [\cdot,\cdot], d_*)$
or, what amounts to the same,
on
$(\Lambda_A D, [\cdot,\cdot]_*, d)$,
are equivalent notions.
This parallels the well known fact that Lie-Rinehart structures on $(A,L)$
are in bijective correspondence with differential graded $R$-algebra
structures on $\roman{Alt}_A(L,A)$.
\smallskip
Let $(A,L',L'')$ be an almost twilled Lie-Rinehart algebra,
having $L'$ and $L''$ finitely generated and projective
as $A$-modules.
The $(A,L')$-module structure on $L''$
induces an
$(A,L')$-module on the dual ${L''}^*$ which, in turn, 
${L''}^*$ being viewed as an abelian Lie algebra and hence
abelian $(R,A)$-Lie algebra, gives rise to
the semi direct product $(R,A)$-Lie algebra
$L'\ltime {L''}^*$.
Likewise the $(A,L'')$-module structure on $L'$
determines the corresponding
semi direct product $(R,A)$-Lie algebra
$L''\ltime {L'}^*$.
Plainly $L=L'\ltime {L''}^*$ and
$D=L''\ltime {L'}^*$ are in duality.
Consider the obvious adjointness isomorphisms
$$
\roman{Alt}_A(L'',\Lambda_A L')
@>>>
\roman{Alt}_A(L''\ltime {L'}^*, A )
=
\roman{Alt}_A(D, A )
\tag4.1.1
$$
and
$$
\Lambda_A L =
\Lambda_A (L'\ltime {L''}^*) 
@>>>
\roman{Alt}_A(L'',\Lambda_A L')
\tag4.1.2
$$
of bigraded $A$-algebras;
these isomorphisms are independent of the 
Lie-Rinehart
semi direct product
constructions and 
instead of $L'\ltime {L''}^*$ and
$L''\ltime {L'}^*$,
we could as well have written
$L'\oplus {L''}^*$
and
$L''\oplus {L'}^*$, respectively.
However, incorporating these semi direct product structures,
we see that, under (4.1.1), 
the Lie-Rinehart differential $d''$
on 
$\roman{Alt}_A(L'',\Lambda_A L')$
passes to the Lie-Rinehart differential
$d_*$ on
$\roman{Alt}_A(D,A)$
and that
under
(4.1.2) 
the (bigraded) Gerstenhaber bracket
$[\cdot,\cdot]$
on
$\Lambda_A L$
passes to the
bigraded Gerstenhaber bracket
$[\cdot,\cdot]'$
$\roman{Alt}_A(L'',\Lambda_A L')$.
Moreover, by construction, the differentials
on both sides of
(4.1.1) are derivations with respect to the multiplicative
structures.

\proclaim{Theorem 4.1}
For an almost twilled Lie-Rinehart algebra
$(A,L',L'')$ having $L'$ and $L''$ finitely generated and projective
as $A$-modules,
$(\roman{Alt}_A(L'',\Lambda_A L'),[\cdot,\cdot]',\parti'')$
is a differential bigraded
Gerstenhaber algebra
if and only if
$(A,L,D)$
is a Lie-Rinehart bialgebra.
\endproclaim

\demo{Proof}
In fact, the first property spelled out above
characterizing  $(A,L,D)$ to be a Lie-Rinehart bialgebra 
is plainly equivalent to 
$(\roman{Alt}_A(L'',\Lambda_A L'),[\cdot,\cdot]',\parti'')$
being a  
differential bigraded
Gerstenhaber algebra. \qed
\enddemo

The following is now immediate, cf. Corollary 4.9 in \cite\twilled.

\proclaim{Corollary 4.2}
An almost twilled Lie-Rinehart algebra
$(A,L',L'')$ having $L'$ and $L''$ finitely generated and projective
as $A$-modules
is a true
twilled Lie-Rinehart algebra
if and only if
$(A,L,D) = (A,L'\ltime {L''}^*,L''\ltime {L'}^*)$
is a Lie-Rinehart bialgebra. \qed
\endproclaim

This result may be proved directly, i.~e. without 
the intermediate
differential bigraded
Gerstenhaber algebra
in (4.1). The reasoning is formally the same, though.
For the special case
where
$L'$ and $L''$ arise from Lie algebroids,
the statement of Corollary 4.2 may be deduced from what is said in
\cite\mackfift.
\smallskip\noindent
{\smc Remark 4.3.} When
$A$ is a field and $\fra g$ an ordinary (finite dimensional)
Lie algebra,
Corollary 4.2 comes down to
the statement that,
in the termiology of
\cite\mackfift, \cite\majidtwo, \cite\mokritwo,
$(\fra g',\fra g'')$
(with the requisite additional structure)
constitutes a matched pair of Lie algebras
(which now amounts to 
$(\fra g',\fra g'')$ being a Lie bialgebra)             
if and only if,
with the obvious structure,
$(\fra g'\ltime {\fra g''}^*,\fra g''\ltime {\fra g'}^*)$
is a
Lie bialgebra.
This fact was known to S. Zakrzewski
\cite\zakrztwo.
It has been spelled out explicitly as Proposition 1 in \cite\stachone.

\medskip

\widestnumber\key{999}
\centerline{References}
\smallskip\noindent

\ref \no \barakont
\by S. Barannikov and M. Kontsevich
\paper Frobenius manifolds and formality of Lie algebras of polyvector fields
\jour Internat. Res. Notices
\vol 4
\yr 1998
\pages 201--215
\finalinfo {\tt alg-geom/9710032}
\endref

\ref \no \batviltw
\by I. A. Batalin and G. S. Vilkovisky
\paper Quantization of gauge theories
with linearly dependent generators
\jour  Phys. Rev. 
\vol D 28
\yr 1983
\pages  2567--2582
\endref

\ref \no \batvilfo
\by I. A. Batalin and G. S. Vilkovisky
\paper Closure of the gauge algebra, generalized Lie equations
and Feynman rules
\jour  Nucl. Phys. B
\vol 234
\yr 1984
\pages  106-124
\endref

\ref \no \batavilk
\by I. A. Batalin and G. S. Vilkovisky
\paper Existence theorem for gauge algebra
\jour Jour. Math. Phys.
\vol 26
\yr 1985
\pages  172--184
\endref

\ref \no \canhawei
\by A. Cannas de Silva, K. Hartshorn, A. Weinstein
\paper Lectures on Geometric Models for Noncommutative Algebras
\paperinfo U of California at Berkeley, June 15, 1998
\endref

\ref \no \cheveile
\by C. Chevalley and S. Eilenberg
\paper Cohomology theory of Lie groups and Lie algebras
\jour  Trans. Amer. Math. Soc.
\vol 63
\yr 1948
\pages 85--124
\endref

\ref \no \gersthtw
\by M. Gerstenhaber
\paper The cohomology structure of an associative ring
\jour Ann. of Math.
\vol 78
\yr 1963
\pages  267-288
\endref

\ref \no \geschthr
\by M. Gerstenhaber and S. D. Schack
\paper Algebras, bialgebras, quantum groups and algebraic
deformations
\paperinfo In: Deformation theory and quantum groups with
applications to mathematical physics, M. Gerstenhaber and J. Stasheff, eds.
\jour Cont. Math.
\vol 134
\pages 51--92
\publ AMS
\yr 1992
\publaddr Providence 
\endref

\ref \no \getzltwo
\by E. Getzler
\paper Batalin-Vilkovisky algebras and two-dimensional topological field
theories
\jour Comm. in Math. Phys.
\vol 195
\yr 1994
\pages 265--285
\endref

\ref \no \herzone
\by J. Herz
\paper Pseudo-alg\`ebres de Lie
\jour C. R. Acad. Sci. Paris 
\vol 236
\yr 1953
\pages 1935--1937
\endref

\ref \no \poiscoho
\by J. Huebschmann
\paper Poisson cohomology and quantization
\jour J. f\"ur die Reine und Angew. Math.
\vol 408
\yr 1990
\pages 57--113
\endref

\ref \no \duality
\by J. Huebschmann
\paper 
Duality for Lie-Rinehart algebras and the modular class
\finalinfo {\tt dg-ga/9702008}
\jour Journal reine angew. Math. 
\vol 510
\yr 1999
\pages 103--159
\endref

\ref \no \bv
\by J. Huebschmann
\paper Lie-Rinehart algebras, Gerstenhaber algebras, and Batalin-
Vilkovisky algebras
\jour Annales de l'Institut Fourier
\vol 48
\yr 1998
\pages 425--440
\endref

\ref \no \extensta
\by J. Huebschmann
\paper 
Extensions of Lie-Rinehart algebras and the Chern-Weil construction
\paperinfo in: Festschrift in honor of J. Stasheff's 60th birthday
\jour Cont. Math. 
\vol 227
\yr 1999
\pages 145--176
\publ Amer. Math. Soc.
\publaddr Providence R. I.
\finalinfo {\tt dg-ga/9706002}
\endref

\ref \no \crosspro
\by J. Huebschmann
\paper Crossed products and twilled Lie-Rinehart algebras
\paperinfo in preparation
\endref

\ref \no \twilled
\by J. Huebschmann
\paper Twilled Lie-Rinehart algebras and differential Batalin-Vilkovisky 
algebras
\paperinfo {\tt math.DG/9811069}
\endref

\ref \no \liribi
\by J. Huebschmann
\paper The modular class and master equation for Lie-Rinehart bialgebras
\paperinfo in preparation
\endref

\ref \no \huebstas
\by J. Huebschmann and J. D. Stasheff
\paper Formal solution of the master equation via HPT and
deformation theory
\jour Forum mathematicum 
\vol 14 
\yr 2002 
\pages 847--868
\finalinfo {\tt math.AG/9906036}
\endref

\ref \no \kosmathr
\by Y. Kosmann-Schwarzbach 
\paper Exact Gerstenhaber algebras and Lie bialgebroids
\jour  Acta Applicandae Mathematicae
\vol 41
\yr 1995
\pages 153--165
\endref

\ref \no \kosmafou
\by Y. Kosmann-Schwarzbach 
\paper From Poisson algebras to Gerstenhaber algebras
\jour Annales de l'Institut Fourier
\vol 46
\yr 1996
\pages 1243--1274
\endref

\ref \no \kosmafiv
\by Y. Kosmann-Schwarzbach 
\paper The Lie bialgebroid of a Poisson-Nijenhius manifold
\jour  Letters in Math. Physics
\vol 38
\yr 1996
\pages 421--428
\endref

\ref \no \kosmagtw
\by Y. Kosmann-Schwarzbach and F. Magri 
\paper Poisson-Lie groups and complete integrability. I.
Drinfeld bigebras, dual extensions and their
canonical representations
\jour  Annales Inst. H. Poincar\'e S\'erie A (Physique th\'eorique)
\vol 49
\yr 1988
\pages 433--460
\endref

\ref \no \koszulon
\by J. L. Koszul
\paper Crochet de Schouten-Nijenhuis et cohomologie
\jour Ast\'erisque,
\vol hors-s\'erie,
\yr 1985
\pages 251--271
\paperinfo in E. Cartan et les Math\'ematiciens d'aujourd'hui, 
Lyon, 25--29 Juin, 1984
\endref

\ref \no \liazutwo
\by B. H. Lian and G. J. Zuckerman
\paper New perspectives on the BRST-algebraic structure
of string theory
\jour Comm. in Math. Phys.
\vol 154
\yr 1993
\pages  613--646
\endref

\ref \no \luweinst
\by J.-H. Lu and A. Weinstein
\paper Poisson Lie groups, dressing transformations, and Bruhat decompositions
\jour J. of Diff. Geom.
\vol 31
\yr 1990
\pages 501--526
\endref

\ref \no \mackfift
\by K. Mackenzie
\paper Double Lie algebroids and the double of a Lie bialgebroid
\paperinfo {\tt math.DG/9808081}
\endref

\ref \no \mackxu
\by K. C. Mackenzie and P.  Xu
\paper Lie bialgebroids and Poisson groupoids
\jour Duke Math. J.
\vol 73
\yr 1994
\pages 415--452
\endref

\ref \no \maclaboo
\by S. Mac Lane
\book Homology
\bookinfo Die Grundlehren der mathematischen Wissenschaften
 No. 114
\publ Springer
\publaddr Berlin $\cdot$ G\"ottingen $\cdot$ Heidelberg
\yr 1963
\endref

\ref \no \majidtwo
\by S. Majid
\paper Matched pairs of Lie groups associated to solutions
of the Yang-Baxter equation
\jour Pac. J. of Math.
\vol 141
\yr 1990
\pages 311--332
\endref

\ref \no \maninfiv
\by Yu. I. Manin
\paper Three constructions of Frobenius manifolds
\paperinfo {\tt math.QA/9801006}
\endref

\ref \no \mokritwo
\by T. Mokri
\paper Matched pairs of Lie algebroids
\jour Glasgow Math. J.
\vol 39
\yr 1997
\pages 167--181
\endref

\ref \no \palaione
\by R. S. Palais
\paper The cohomology of Lie rings
\jour  Proc. Symp. Pure Math.
\vol III
\yr 1961
\pages 130--137
\paperinfo Amer. Math. Soc., Providence, R. I.
\endref

\ref \no \rinehone
\by G. Rinehart
\paper Differential forms for general commutative algebras
\jour  Trans. Amer. Math. Soc.
\vol 108
\yr 1963
\pages 195--222
\endref

\ref \no \schecone
\by V. Schechtman
\paper Remarks on formal deformations and Batalin-Vilkovisky algebras
\paperinfo {\tt math.AG/9802006}
\endref

\ref \no \stachone
\by P. Stachura
\paper Double Lie algebras and Manin triples
\paperinfo {\tt q-alg/9712040}
\endref

\ref \no \stashnin
\by J. D. Stasheff
\paper Deformation theory and the Batalin-Vilkovisky 
master equation
\paperinfo in: Deformation Theory and Symplectic Geometry,
Proceedings of the Ascona meeting, June 1996,
D. Sternheimer, J. Rawnsley, S.  Gutt, eds.,
Mathematical Physics Studies, Vol. 20
\publ Kluwer Academic Publishers
\publaddr Dordrecht/Boston/London
\yr 1997
\pages 271--284
\endref

\ref \no \tianone
\by G. Tian
\paper A note on Kaehler manifolds with $c_1=0$
\paperinfo preprint
\endref

\ref \no \todortwo
\by A. N. Todorov
\paper The Weil-Petersson geometry
of the moduli space of $\fra {su}(n \geq 3)$ (Calabi-Yau)
manifolds, I.
\jour Comm. Math. Phys.
\vol 126
\yr 1989
\pages 325--346
\endref

\ref \no \wittetwe
\by E. Witten
\paper Mirror manifolds and topological field theory
\paperinfo in: Essays on mirror manifolds, 
S. T. Yau, ed.
\publ International Press Co.
\publaddr Hong Kong
\yr 1992
\pages 230--310
\endref

\ref \no \xuone
\by P. Xu
\paper 
Gerstenhaber algebras and BV-algebras
in Poisson geometry
\paperinfo preprint, 1997
\endref

\ref \no \zakrztwo
\by S. Zakrzewski
\paper Poisson structures on the Poincar\'e groups
\jour Comm. Math. Phys. 
\vol 185
\yr 1997
\pages 285--311
\endref

\enddocument